\documentclass{article} 
\usepackage{texdraw}\usepackage{graphicx}
\usepackage{amsmath} 
\usepackage{color}
\usepackage{colordvi}
\usepackage{texdraw}\usepackage{graphicx}
\usepackage{dcpic,pictex}

\definecolor{purple}{cmyk}{0.45,0.86,0,0}\definecolor{brickred}{cmyk}{0,0.89,0.94,0.28}
\definecolor{maroon}{cmyk}{0,0.87,0.68,0.32}\definecolor{lyellow}{cmyk}{0,0,0.68,0}
\definecolor{magenta}{cmyk}{0,1,0,0}
\definecolor{zals}{cmyk}{0.75,0.12,0.98,0.5}
\definecolor{pamats}{cmyk}{0.54,0.00,0.84,0.0}
\definecolor{brickred}{cmyk}{0,0.89,0.94,0.28}
\definecolor{cadetblue}{cmyk}{0.62,0.57,0.23,0}
\definecolor{teo}{cmyk}{0.63,0.83,0.33,0.13}

\def\newpic#1{%   
   \def\emline##1##2##3##4##5##6{%
      \put(##1,##2){\special{em:point #1##3}}% 
      \put(##4,##5){\special{em:point #1##6}}%
      \special{em:line #1##3,#1##6}}}
\newpic{}

\usepackage[cp1251,cp1257]{inputenc}\usepackage[english,latvian]{babel}\usepackage{amssymb,latexsym,amsfonts}
\newcommand{\apz}{\protect\makebox[1.8\width]{\rule[3.3pt]{8.5pt}{0.4pt} \hspace*{-12.6pt}\protect\raisebox{-.3\height}{$\leftharpoondown$}}}

\newtheorem{teor}{Theorem}[section]

\newtheorem{corollary}[teor]{Corollary}
\newtheorem{lemma}[teor]{Lemma}
\newtheorem{proposition}[teor]{Proposition}
\newtheorem{definition}[teor]{Definition}

\newtheorem{example}[teor]{Example}

\input cyracc.def

\begin{document}
\selectlanguage{english}
%1-----------------------------------------------------------------------------------
\title{The infinite dihedral group}
\date{}
\author{J\=anis Buls \\
{\small Department  of  Mathematics, University of Latvia, Jelgavas iela 3,}\\
{\small R\=\i ga, LV-1004 Latvia,
buls@edu.lu.lv
}}
\def\keywords{\begin{center}{\bf Keywords}\end{center}
{automata (machines) groups, the infinite dihedral group}}

\maketitle
\title{}
% ************* Begin of Abstract ******************%
\abstract{\small{We describe the infinite dihedral group as automaton group. We collect basic results and give full proofs in details for all statements. }}
% ************* End of Abstract ******************%


\section{ Preliminaries} 
%%%%%%%%%%%%%%%%%%%%%%%%%%%%%%%%%%%%%%%%%%%%%%%%%%%%%%%%%%%%%%%%%%%%%%%%%

We use standard conventions.
Nethertheless we recall some notation. For more details see [1].

Let $A$ be a finite non-empty set and $A^*$ be
the free monoid generated by $A$. The set $A$ is also called an {\em 
alphabet}, its  elements are called {\em letters} and 
those of $A^*$ are called {\em finite words}. The identity element of $A^*$ is called 
an {\em empty word} and denoted by $\lambda$. We set 
$A^+=A^*\backslash\{\lambda\}$.

\begin{definition}
A 2-sorted algebra $\mathcal{M}=\langle Q,A,\circ ,\ast \rangle$ is called  a 
{\em Mealy machine} if $Q,A$ are finite, nonempty sets, the mappings 
$Q\times A \stackrel{\circ}{\longrightarrow}Q$,
$Q\times A \stackrel{\ast }{\longrightarrow}A$
are  total  functions. 
\end{definition}

\begin{itemize}
	\item $Q$ is called  a set of states;
	\item $A$  --- an alphabet;
	\item $\circ$ --- the transition function;
	\item $*$ --- the output function.
\end{itemize}
 One draws 
\begin{eqnarray}
\bigcirc &\xrightarrow{a/b}& \bigcirc \label{for1}\\
q_1 & & q_2 \nonumber
\end{eqnarray}
%\[\underset{\displaystyle q_1}{\bigcirc} \xrightarrow{a/b} \underset{\displaystyle q_2}{\bigcirc}\]
\medskip\\
to mean $q_1 \circ a= q_2, q_1 *a = b$. That is, in state $q_1$ with input $a$, the machine
outputs $b$ and goes to state $q_2$.

The machine is said to be {\em invertible} if for each state $q\in Q$ the output function is a permutation.
The inverse machine is the
machine $\mathcal{M}^{-1}$ obtained by switching the input $a$ and the output $b$ on each
arrow (\ref{for1}). The resulting transition in the inverse automaton is
(\ref{for2}). 

\begin{eqnarray}
\bigcirc &\xrightarrow{b/a}& \bigcirc \label{for2}\\
q_1 & & q_2 \nonumber
\end{eqnarray}

The mappings $\circ$ and $\ast$ may be extended to $Q\times A^*$  by defining
\[
\begin{array}{lr}
q\circ \lambda =q,\quad & q\circ (ua)=(q\circ u)\circ a, \\
q\ast \lambda  =\lambda ,\quad & q\ast (ua)=(q\ast u)\#((q\circ u)\ast a)\,,
\end{array}
\] 
for each $q\in Q$, $(u,a)\in A^*\times A$. Here $\#$ means the concatenation of words.

Let $\mathbb{N}=\{0,1,2,\ldots, n, \ldots\}$.
An (indexed) infinite word $x$ on the alphabet $A$ is any total mapping 
$x\,:\,\mathbb{N}\rightarrow A$. We shall set for any $i\ge0$, $x_i=x(i)$
and write
\[
x=(x_i)=x_0x_1\cdots x_n\cdots \;.
\]
The set of all the infinite words over $A$ is denoted by $A^\omega$. 

Let $A^\infty=A^*\cup A^\omega$.
The mapping $\ast$ may be extended to $Q\times A^\infty$  by defining
\[
q*x=\lim\limits_{n\to \infty}q*x[0,n],
\]
where $x[0,n]=x_0x_1\cdots x_n$.

We like to be more flexible therefore we sometimes write $uq$ instead of $q*u$ and the state $q$ identify with the map $q*u$. In other words we identify the map $\bar q: A^\infty \to A^\infty: u\mapsto q*u$ with $q$. 

The semigroup  generated by the machine $\mathcal{M}$, which we shall call the machine semigroup of $\mathcal{M}$, is the semigroup $\langle \mathcal{M} \rangle_+$
generated by $\{ \bar q\,|\, q\in Q\}$. If $\mathcal{M}$ is invertible, then the group generated by the machine $\mathcal{M}$, which we shall call the machine group of $\mathcal{M}$, is the group $\Gamma(\mathcal{M})$ generated by $\{ \bar q\,|\, q\in Q\}$.

\section{The infinite dihedral group}
%%%%%%%%%%%%%%%%%%%%%%%%%%%%%%%%%%%%%%%%%%%%%%%%%%%%%%%%%%%%

Let $\mathbb{Z}$ be  the set of integers and $\mathfrak{S}(\mathbb{Z})$ be the symmetric group on $\mathbb{Z}$.
\begin{definition}
\[
D_\infty=\{s\in\mathfrak{S}(\mathbb{Z})\,|\,|s(i)-s(j)|=|i-j|\}
\]
\end{definition}

\begin{proposition}
$D_\infty\le\mathfrak{S}(\mathbb{Z})$
\end{proposition}

$\Box$ (i) Let $s,\sigma\in D_\infty$, then
\[
|s(\sigma(i))-s(\sigma(j))|=|\sigma(i)-\sigma(j)|=|i-j|
\]
Hence the composition of bijections $\sigma s\in D_\infty$.

(ii) Since $s$ is a bijection then exist  $k$ and $n$ such that $s(k)=i$ and $s(n)=j$. Hence
\[
|s^{-1}(i)-s^{-1}(j)|=|s^{-1}(s(k))-s^{-1}(s(n))|=|k-n|=|s(k)-s(n)|=|i-j|
\]
Therefore  $s^{-1}\in D_\infty$.

This completes the proof that $D_\infty$ is a subgroup of $\mathfrak{S}(\mathbb{Z})$.
\rule{2mm}{2mm}

\begin{definition}
$D_\infty$ is called the infinite dihedral group.
\end{definition}

\begin{lemma}\label{l5.19.10}
$s\in D_\infty \Leftrightarrow \exists n\forall i \;[ s(i)=i+n \vee s(i)=-i+n]$
\end{lemma}

$\Box \Leftarrow \;$ (i)  $|s(i)-s(j)|=|i+n-(j+n)|=|i-j|$. Hence $s\in D_\infty$.

(ii) $|s(i)-s(j)|=|-i+n-(-j+n)|=|j-i|=|i-j|$. Therefore $s\in D_\infty$.

$\Rightarrow $ Let $s\in D_\infty$ then $|s(i+1)-s(i)|=|i+1-1|=1$.  Hence $s(i+1)=s(i)+1$ or $s(i+1)=s(i)-1$.

(i) Let $s(i+1)=s(i)+1$ but $s(i+2)=s(i+1)-1$ then 
\[s(i+2)=s(i+1)-1=s(i)+1-1=s(i).\]
 Therefore               $|s(i+2)-s(i)|=0\ne |i+2-i|=2$.
 Corollary: if $s(i+1)=s(i)+1$ then $\forall j\>\ge i\; s(j+1)=s(j)+1$.

(ii) Let $s(i+1)=s(i)-1$ but $s(i+2)=s(i+1)+1$ then 

\[
s(i+2)=s(i+1)+1=s(i)-1+1=s(i). 
\]
Therefore $|s(i+2)-s(i)|=0\ne |i+2-i|=2$.
 Corollary: if $s(i+1)=s(i)-1$ then $\forall j\>\ge i\; s(j+1)=s(j)-1$.

Hence (from (i) and (ii))  $\forall j\; s(j)=s(j-1)+1$, or $\forall j\; s(j-1)=s(j)-1$.
In the first case, there exists $n$ such that $\forall i\; s(i)=i+n$. In the second case, there exists $n$ such that
$\forall i\; s(i)=-i+n$. Note $n=s(0)$.
\rule{2mm}{2mm}

\begin{corollary}
$D_\infty$ is not a commutative group.
\end{corollary}

$\Box$  Let $s:i\mapsto i+1 ,\sigma: i\mapsto -i+2$, then 
\begin{eqnarray*}
s(\sigma(1)) &=& s(-1+2)=s(1)=2,\\
\sigma(s(1)) &=& \sigma (2)=0. \quad  \rule{2mm}{2mm}
\end{eqnarray*} 

\begin{lemma}
The mapping $\varphi: \mathbb{Z}_2\to \mathfrak{A}ut(\mathbb{Z})$, where 
\begin{eqnarray*}
\varphi^0 &=& \mathbb{I},\\
\varphi^1 &=& -\mathbb{I},
\end{eqnarray*}
and  $\mathbb{I}:\mathbb{Z}\to \mathbb{Z}: x\mapsto x$ is an identity map, 
is a  group  homomorphism.
\end{lemma}

$\Box$ We specify that the group $\mathbb{Z}_2$ is meant to be the cyclic group of order 2 $\langle \mathbb{Z}_2, +\rangle$ but the  group $\mathbb{Z}$ is meant to be the group of integers under addition. 

We write the automorphism on $\mathbb{Z}$ associated to $k$ as $\varphi^k$,
so $\varphi^k : \mathbb{Z} \to \mathbb{Z}$ is a bijection and $\varphi^k(hh')=\varphi^k(h)\circ \varphi^k(h')$
 for all $h,h'$ in $\mathbb{Z}$.

(i)  $\varphi^{0+0}=\varphi^0=\mathbb{I}=\mathbb{I}\circ \mathbb{I}=\varphi^0\circ\varphi^0$.

(ii) $\varphi^{0+1}=\varphi^1=-\mathbb{I}=\mathbb{I}\circ(- \mathbb{I})=\varphi^0\circ\varphi^1$.

(iii)  $\varphi^{1+0}=\varphi^1=-\mathbb{I}=-\mathbb{I}\circ \mathbb{I}=\varphi^1\circ\varphi^0$.

(iv)  $\varphi^{1+1}=\varphi^0=\mathbb{I}=-\mathbb{I}\circ(- \mathbb{I})=\varphi^1\circ\varphi^1$.
\rule{2mm}{2mm}

Let 
\begin{itemize}
	\item[$\bullet$] $\mathcal{N}$ and $\mathcal{H}$ be groups;
	\item[$\bullet$] $\varphi: \mathcal{H} \to \mathfrak{A}ut(\mathcal{N})$ ---  a group homomorphism, where 
	$\mathfrak{A}ut(\mathcal{N})$ --- the automorphism group of the group $\mathcal{N}$.
\end{itemize}
We shall write $\varphi^h\in\mathfrak{A}ut(\mathcal{N})$ insted of $\varphi(h)$.
Then define $\mathcal{G}=\mathcal{N}\rtimes_\varphi \mathcal{H}=\mathcal{N}\rtimes \mathcal{H}$ to be the set
$\mathcal{N}\times \mathcal{H}$ with the multiplication defined by
\[
(n_1,h_1)(n_2,h_2)=(n_1\varphi^{h_1}(n_2),h_1h_2).
\]

\begin{definition}\label{def9.1}
Let $\mathcal{G}$ be a group with subgroups $\mathcal{N}$ and $\mathcal{H}$ such that 
\[
\mathcal{N}\cap \mathcal{H}=\{e\} \quad {\rm and} \quad \mathcal{N}\mathcal{H}=\mathcal{G}.
\]
If $\mathcal{N}$ is normal (but not necessarily $\mathcal{H}$), then we say that $\mathcal{G}$ is the semidirect product of $\mathcal{N}$ and $\mathcal{H}$. Here $e$ is the neutral element of $\mathcal{G}$.
\end{definition}

Other necessary details see for example in \cite{lamp}.

\begin{proposition}
$\mathbb{Z} \rtimes \mathbb{Z}_2\cong D_\infty$
\end{proposition}

$\Box$ 
Let $(k,\alpha), (n,\beta)\in \mathbb{Z} \rtimes \mathbb{Z}_2$ then from previous Lemma 
\[(k,\alpha)(n,\beta)=(k+\varphi^\alpha(n), \alpha+\beta).\]

Choose
$\psi:\mathbb{Z} \rtimes \mathbb{Z}_2\to D_\infty: (k,\alpha)\mapsto s$ where
\[
s:\mathbb{Z}\to \mathbb{Z}: i\mapsto \begin{cases}
i+k, & {\rm ja}\; \alpha=0;\\
-i+k, & {\rm ja}\; \alpha=-1. 
\end{cases} 
\]
 Now by Lemmma \ref{l5.19.10} we conclude that  $\psi$ is  a bijection. Its remains to proof that  $\psi $ is a group homomorphism.

Let  $\psi(k,\alpha)=s$ and $\psi(n,\beta)=\sigma$.

(i) If $\alpha=\beta=0$ then  $s:i\mapsto i+k$ and $\sigma: i\mapsto i+n$. Hence $(k,\alpha)(n,\beta)=(k+n,\alpha+\beta)=(k+n,0)$ and $\psi((k,\alpha)(n,\beta)): i\mapsto i+k+n$.
Otherwise $s\circ \sigma:i\mapsto s(i+n)=i+n+k$. Thence
\centerline{ $\psi((k,\alpha)(n,\beta))=\psi(k,\alpha)\circ\psi(n,\beta)$.}

(ii) If $\alpha=0$ but $\beta=1$ then  $s:i\mapsto i+k$ and $\sigma: i\mapsto -i+n$. Hence $(k,\alpha)(n,\beta)=(k+n,\alpha+\beta)=(k+n,1)$ and $\psi((k,\alpha)(n,\beta)): i\mapsto -i+k+n$.
Otherwise $s\circ \sigma:i\mapsto s(-i+n)=-i+n+k$. Thence \\
 \centerline{$\psi((k,\alpha)(n,\beta))=\psi(k,\alpha)\circ\psi(n,\beta)$.}

(iii) If $\alpha=1$ but $\beta=0$ then  $s:i\mapsto -i+k$ and $\sigma: i\mapsto i+n$. Hence $(k,\alpha)(n,\beta)=(k-n,\alpha+\beta)=(k-n,1)$ and $\psi((k,\alpha)(n,\beta)): i\mapsto -i+k-n$.
Otherwise $s\circ \sigma:i\mapsto s(i+n)=-i-n+k$. Thence \\
\centerline{$\psi((k,\alpha)(n,\beta))=\psi(k,\alpha)\circ\psi(n,\beta)$.}

(iv) If $\alpha=\beta=1$ then  $s:i\mapsto -i+k$ and $\sigma: i\mapsto -i+n$. Hence $(k,\alpha)(n,\beta)=(k-n,\alpha+\beta)=(k-n,0)$ and $\psi((k,\alpha)(n,\beta)): i\mapsto i+k-n$.
Otherwise $s\circ \sigma:i\mapsto s(-i+n)=i-n+k$. Thence \\
\centerline{$\psi((k,\alpha)(n,\beta))=\psi(k,\alpha)\circ\psi(n,\beta)$.}

Summing up (cases (i)--(iv)),  $\psi$ is a group homomorphism.
\rule{2mm}{2mm}

\begin{lemma}
Let $e$ be the identity of  group $\mathbb{Z} \rtimes \mathbb{Z}_2$ then $x^2=e$ and $xax=a^{-1}$ where $x=(-1,1)$ and $a=(1,0)$.
\end{lemma}

$\Box$ $x^2=(-1,1)(-1,1)=(-1+\varphi^1(-1),1+1)=(-1+1,0)=(0,0)=e$
\begin{eqnarray*}
xax &=& (-1,1)(1,0)(-1,1)=(-1+\varphi^1(1),1+0)(-1,1)\\
&=& (-1-1,1)(-1,1)=(-2,1)(-1,1)\\
&=& (-2+\varphi^1(-1),1+1)=(-2+1,0)=(-1,0)
\end{eqnarray*}

Let us make sure that $(-1,0)=a^{-1}$.

$(1,0)(-1,0)=(1+\varphi^0(-1),0+0)=(1-1,0)=(0,0)=e$ \rule{2mm}{2mm}

\begin{proposition}
$\mathbb{Z} \rtimes \mathbb{Z}_2=\langle x, a\rangle$, i.e., the set $\{x, a\}$ generates \\
$\mathbb{Z} \rtimes \mathbb{Z}_2$.
\end{proposition}

$\Box$ (i) $a^2=(1,0)(1,0)=(1+\varphi^0(1),0+0)=(1+1,0)=(2,0)$. 

The rest is induction. Let $a^n=(n,0)$. Hence
\[
a^{n+1}=(a^na)=(n,0)(1,0)=(n+\varphi^0(1),0+0)=(n+1,0).
\]

(ii) $a^n(-n,0)=(n,0)(-n,0)=(n+\varphi^0(-n),0+0)=(n-n,0)=(0,0)=e$. Since $a^na^{-n}=e$ then $a^{-n}=(-n,0)$.

Therefore $\forall n\in \mathbb{Z} \; (n,0)\in \langle x, a\rangle $.

(iii) $ax=(1,0)(-1,1)=(1+\varphi^0(-1),0+1)=(1-1,0)=(0,1)$.
  \[
a^nx=(n,0)(-1,1)=(n+\varphi^0(-1),0+1)=(n-1,1)
\]
for all $n\in\mathbb{Z}$. Thence $\forall k\in \mathbb{Z} \; (k,1)\in \langle x, a\rangle $.

We have demonstrated: if $(n,\alpha) \in \mathbb{Z} \rtimes \mathbb{Z}_2$ then $(n,\alpha)\in  \langle x, a\rangle $.
\rule{2mm}{2mm}

\begin{lemma}\label{l5.19.16}
Let  $a\ne e\ne x$ and $x\ne a$ are elements of a group\;  $\langle x, a \rangle$, furthermore $x^2=e$ and $xax=a^{-1}$ where  $e$ is the identety of  $\langle x, a \rangle$.
 If  $\langle x, a \rangle$ is an infinite group then $\mathbb{Z} \rtimes \mathbb{Z}_2\cong\langle x, a\rangle$.
\end{lemma}

$\Box$ 
 (i) Assume first that $a^s=e$ for natural $s$. Hence we have two sets 

\centerline{$A_0\apz\{e,a, a^2,\ldots, a^{s-1}\}$ and $X_0\apz\{x,ax, a^2x,\ldots, a^{s-1}x\}$.}

\begin{itemize}
\item If $u=a^k$ and $v=a^n$ then $uv=a^{k+n}\in A_0$.
\item If $u=a^k$ and $v=a^nx$ then $uv=a^k a^nx=a^{k+n}x\in X_0$.
\item If $u=a^kx$ and $v=a^n$ then 
\[uv=a^kxa^n=a^kxa^n xx=a^ka^{-n}x=a^{k-n}x\in X_0.\]

The identity $xa^nx=a^{-n}$ can be proven:
\[
xa^nx=xaxxa^{n-1}x=a^{-1}a^{-n+1}=a^{-n}
\]
\item If $u=a^kx$ and $v=a^nx$, then $uv=a^kxa^nx=a^ka^{-n}=a^{k-n}\in  A_0$.
\end{itemize}
Hence $|\langle x, a \rangle|\le|A_0\cup X_0|\le 2s$. Contradiction!

(ii) $X\apz\{ax^k \,|\, k\in \mathbb{Z}\}$ is an infinite set.

Proof by contradiction. Assume $\exists kn\; a^kx=a^nx$ then by the  cancellation law  $a^{k-n}=x$, and therefore   $a^{2(k-n)}=x^2=e$. This is contradistion (we know that $A\apz\{a^k\,|\, k\in\mathbb{Z}\}$ is an infinite set).

(iii) $A\cap X=\emptyset$.

Proof by contradiction. Assume $\exists kn\; a^kx=a^n$ then by the  cancellation law  $a^{n-k}=x$, and therefore   $a^{2(k-n)}=x^2=e$. This is contradistion (we know that $A$ is an infinite set).

(iv) $G=A\cup X$. 

Since $G=\langle x, a\rangle$ then $G\supseteq A\cup X$.
We shall copy further decision from case (i). 
 Let $u,v \in A\cup X$.

\begin{itemize}
\item If $u=a^k$ and $v=a^n$ then $uv=a^{k+n}\in A$.
\item If $u=a^k$ and $v=a^nx$ then $uv=a^k a^nx=a^{k+n}x\in X$.
\item If $u=a^kx$ and $v=a^n$ then 
\[uv=a^kxa^n=a^kxa^n xx=a^ka^{-n}x=a^{k-n}x\in X.\]
\item If $u=a^kx$ and $v=a^nx$ then $uv=a^kxa^nx=a^ka^{-n}=a^{k-n}\in  A$.
\end{itemize}

(v) Now from cases (iii) and  (iv) 
\[
\varphi: G\to \mathbb{Z} \rtimes \mathbb{Z}_2: a^k\mapsto (k,0), a^kx \mapsto (k-1,1)
\]
we can conclude $\varphi$ is bijective. Its remains to
proof that $\varphi$ is a group homomorphism, namely,  $\varphi(uv)=\varphi(u)\varphi(v)$ for all $u,v\in G$.
\begin{itemize}
\item If $u=a^k$ and $v=a^n$ then $\varphi(uv)=\varphi(a^{k+n})=(k+n,0)$. From the other hand  
$\varphi(u)\varphi(v)=(k,0)(n,0)=(k+n,0)$.
\item If $u=a^k$ and $v=a^nx$ then  $\varphi(uv)=\varphi(a^{k+n}x)=(k+n-1,1)$. From the other hand
$\varphi(u)\varphi(v)=(k,0)(n-1,1)=(k+n-1,1)$.
\item If $u=a^kx$ and $v=a^n$ then $\varphi(uv)=\varphi(a^{k-n}x)=(k-n-1,1)$. From the other hand
$\varphi(u)\varphi(v)=(k-1,1)(n,0)=(k-1-n,1)$.
\item If $u=a^kx$ and $v=a^nx$ then  $\varphi(uv)=\varphi(a^{k-n})=(k-n,0)$. From the other hand
$\varphi(u)\varphi(v)= (k-1,1)(n-1,1)=(k-1-n+1,1+1)=(k-n,0)$. \quad \rule{2mm}{2mm}
\end{itemize}

\begin{lemma}\label{l5.19.17}
Let  $\alpha\ne e\ne \beta$ and $\alpha\ne \beta$ are elements of a group\;  $G\apz\langle \alpha, \beta \rangle$, furthermore $\alpha^2=e=\beta^2$  where  $e$ is the identety of  $G$.
 If  $G$ is an infinite group then $\mathbb{Z} \rtimes \mathbb{Z}_2\cong\langle \alpha, \beta\rangle$.
\end{lemma}

$\Box$  (i) Let $a\apz \alpha\beta$ and $x\apz \beta$ then $\langle a,x\rangle\subseteq\langle \alpha, \beta \rangle$.

Since $ax=\alpha\beta\beta=\alpha$ then $\langle \alpha, \beta \rangle \subseteq \langle a,x\rangle$.

Hence $\langle \alpha, \beta \rangle = \langle a,x\rangle$.

(ii) If it turns out that $a=x$ then $\alpha=ax=xx=\beta^2=e$. Contradiction! Thus $a\ne x$.

(iii) If it turns out that $a=e$ than $e=a=\alpha\beta$. Therefore $\alpha=\alpha e=\alpha (\alpha\beta)=\alpha^2\beta=e\beta=\beta$. Contradiction because $\alpha\ne \beta$.
Thus $a\ne e$.

(iv) $xax=\beta \alpha\beta \beta=(\beta\alpha)(\beta^2)=(\beta\alpha)e=\beta\alpha$ and $a^{-1}=(\alpha\beta)^{-1}=\beta^{-1}\alpha^{-1}$. Since $\alpha^2=e=\beta^2$ then
$\alpha^{-1}=\alpha$ and $\beta^{-1}=\beta$. Hence $a^{-1}=\beta^{-1}\alpha^{-1}=\beta\alpha =xax$.

(v)   We have checked  (cases (i)--(iv)) all requirements from Lemma \ref{l5.19.16}. Since $\langle \alpha, \beta \rangle = \langle a,x\rangle$ then $\mathbb{Z} \rtimes \mathbb{Z}_2\cong\langle a,x\rangle=\langle \alpha, \beta \rangle$.

Note that according to the given the group $ \langle \alpha, \beta \rangle $ is infinite.
\rule{2mm}{2mm}

\section{The infinite dihedral group as machine group}

\begin{figure}[h]
\input{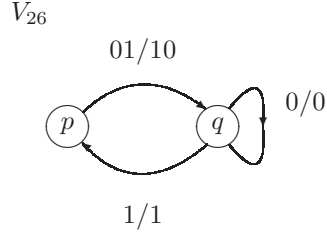}
\caption{Machine defines $D_\infty$.}
\label{z94}
\end{figure}

\begin{example}
$\Gamma(V_{26})\cong D_\infty$
\end{example}

$\Box$  States $pp$ and $qq$ of machine $V'_{26}$ (see fig. \ref{z95}) are equivalent and define maps $\bar p^2$, $\bar q^2$. We have identity $\bar p^2=\bar q^2=\mathbb{I}:\{0,1\}^\infty \to \{0,1\}^\infty: x\mapsto x$. Thus these elements are identity of group $\Gamma(V_{26})$. Its remains  to
proof that $\Gamma(V_{26})$ is an infinite set (see Lemma \ref{l5.19.17}).

\begin{figure}[h]
\input{aut95.pic}
\caption{$V_{26}\leadsto V_{26}$}
\label{z95}
\end{figure}

\begin{figure}[h]
\input{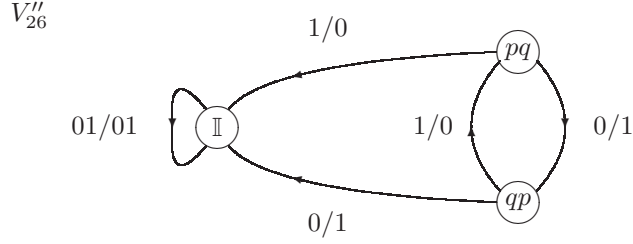}
\caption{The machine's  $V'_{26}$ reduced machine.}
\label{z96}
\end{figure}

\begin{lemma}
Let  $x=x_0x_1\cdots x_n \cdots\in\overline{0,1}^\omega$ then

\centerline{$\forall n\exists s\; x(\overline{pq})^s[n]=\tilde x_n$}
where
\[
\tilde x_n\apz
\begin{cases} 0, &{\rm if}\; x_n=1,\\
1, &{\rm if}\; x_n=0. 
\end{cases}
\]
\end{lemma}

$\Box$ The proof is inductive by the length of the prefix. 

{\bf The induction basis.} We shall base on fig. \ref{z96} 

(i) Let $x=0y$ then  $x\overline{pq}=0y\overline{pq}=1y'$ where $y'\in \overline{0,1}^\omega$.

If $x=1y$ then $x\overline{pq}=1y\overline{pq}=0y$.

If $x=00y$ then $x\overline{pq}=00y\overline{pq}=11y$.

If $x=01y$ then $x\overline{pq}=01y\overline{pq}=10 z$ where $z\in\{0,1\}{0,1}^\omega$.

If $x=10y$ then $x(\overline{pq})^2=10y(\overline{pq})^2=00y\overline{pq}=11y$.

If $x=11y$ then $x(\overline{pq})^2=11y(\overline{pq})^2=01y\overline{pq}=10z$.

This had proved the induction basis.

(ii) We should like to use less  cumbersome notation. Thus  $x\mapsto y$ will mean 
\[
\exists s\; x(\overline{pq})^s=y.
\]

(iii) {\bf Indukcion step}: if $|u|=n$ then $u0x\mapsto u'1x'$ and $u1x\mapsto u''0x''$ where $x', x''$ are $\omega$--words and $|u|=|u'|=|u''|$.

a) Let $|u|=2k-1=n$ --- an odd number, then
\[
u0x\mapsto u'1x' \mapsto v0y,\quad |u|=|u'|=|v| 
\]
where $x'$, $y$ are  $\omega$--words. Assume this is the least s that satisfies
requirement: 1 is replaced by 0, namely, 
\[
u'1x'(\overline{pq})^s=v0y
\]
but $u'1x'(\overline{pq})^{s-1}=v'1y'$ where $|u'|=|v'|$. Therefore
 $v'1y'\overline{pq}=v0y$.
This is possible if 
\[
v'1y'\overline{pq}=v'\overline{pq}\#(pq\circ v')*1\# y.
\]
Since $|v'|=2k-1$ then $pq\circ v'= qp$ or  it is $\mathbb{I}$ (see fig. \ref{z96}). The state $\mathbb{I}$ does not change a letter. Hence $pq\circ v'= qp$. This demonstrates that $pq\circ v'1=pq$. There is a change of a letter in the state $pq$, i.e., the first letter of $y'$ is not the first letter of $y$. Therefore exists $\sigma$ such that
$u0x[n+1]\ne u0x(\overline{pq})^\sigma[n+1]$.

b) Let $|u|=2k=n$ --- an even number, then
\[
u1x\mapsto u''0x'' \mapsto v1y,\quad |u|=|u''|=|v| 
\]
where $x''$, $y$ are  $\omega$--words.  Assume this is the least s that satisfies
requirement: 0 is replaced by 1, namely,  
\[
u''0x''(\overline{pq})^s=v1y
\]
but $u''0x''(\overline{pq})^{s-1}=v'0y'$ where $|u''|=|v'|$. Therefore $v'0y'\overline{pq}=v1y$.
This is possible if 
\[
v'0y'\overline{pq}=v'\overline{pq}\#(pq\circ v')*0\# y.
\]
Since $|v'|=2k$ then $pq\circ v'= pq$ or  it is $\mathbb{I}$ (see fig. \ref{z96}). The state $\mathbb{I}$ does not change a letter. Hence $pq\circ v'= pq$. This demonstrates that $pq\circ v'1=qp$.
There is a change of a letter in the state $qp$, i.e., the first letter of $y'$ is not the first letter of $y$. Therefore exists $\sigma$ such that
$u1x[n+1]\ne u1x(\overline{pq})^\sigma[n+1]$.
\rule{2mm}{2mm}
\medskip

Note (see fig.  \ref{z96} zīm.) $u00\circ pq=\mathbb{I}$. Therefore $u0^\omega\overline{pq}=v0^\omega$ where $|v|\le|u|+2$.
Hence $(\overline{pq})^s0^\omega=w0^\omega$ where $|w|\le 2s$. This means that $(\overline{pq})^s0^\omega[2s]=0$.\par
If it turns out that $(\overline{pq})^s0^\omega=(\overline{pq})^{s+n}0^\omega$ such that $n> 0$ then there would be
a sequence of words 
\[
u_0, u_1,\ldots, u_s, \ldots, u_{s+n}
\]
where $\forall i\in (\overline{0,n-1})(\overline{pq})^{s+n+i}0^\omega=u_{s+i}0^\omega$ besides $u_{s+i}=u_{s+n+i}$.
As a result
\[
\forall k\;  (\overline{pq})^{s+k}0^\omega\in\{u_00^\omega,u_10^\omega, \ldots, u_{s+n-1}0^\omega\}.
\]

We know $(\overline{pq})^{s+i}0^\omega[2(s+i)]=0$. Therefore $\forall i\in \overline{0,s+n-1} \;  |u_i|\le 2(s+n)$. 
We have
\[
\forall k\;  (\overline{pq})^{s+k}\in\{u_00^\omega,u_10^\omega, \ldots, u_{s+n-1}0^\omega\}.
\] 
Hence
$\forall j \; ( \overline{pq})^j[2(s+n+1)]=0 $. Contradiction because the Lemma above tells us: exists $\varkappa$ such that $(\overline{pq})^\varkappa0^\omega[2(s+n+1)]=1$. Therefore, all words $(\overline{pq})^j0^\omega$ are different. This means that all maps $(\overline{pq})^j$ are different. Summing up, $\Gamma(V_{26})$ is an infinite set.
\rule{2mm}{2mm}

\end{document}